# Probability theory and its models

**Paul Humphreys**[1]

*University of Virginia*

**Abstract:** This paper argues for the status of formal probability theory as a mathematical, rather than a scientific, theory. David Freedman and Philip Stark's concept of model based probabilities is examined and is used as a bridge between the formal theory and applications.

## 1. Introduction

David Freedman's work has an unusual breadth, ranging from monographs in mathematical probability through results in theoretical statistics to matters squarely in applied statistics. Over the years, we have had many philosophically flavored discussions about each of these areas and the content of some of these discussions can be distilled into three questions. First, is probability theory a mathematical theory or does it, in virtue of its wide applicability in various areas of science, count as scientific? Secondly, how do we get from the abstract theory of probability to the world? Thirdly, what, if any, are the correct interpretations of probability? Material to answer the second of these questions can be found in Freedman's many publications on statistical models, such as Freedman [7, 8]. Some answers to the third question are directly addressed in Freedman [6]. The answer to the first question may seem obvious to those who work in mathematical probability and it is to Freedman – probability theory is a part of pure mathematics. Ultimately, I think that this is the correct answer but I do not think that it is quite as obvious as it might appear. In what follows I shall try to provide a unified approach to all three questions.[2]

## 2. The status of probability theory

My first question is a little blunt but the founder of modern probability theory addressed it succinctly. For Kolmogorov, probability theory was to be a part of mathematics: "The author set himself the task of putting in their natural place, among the general notions of modern mathematics, the basic concepts of probability theory..." (Kolmogorov [15], p. v).[3] It is sometimes suggested that in constructing the measure theoretic account of probability, Kolmogorov solved one half of

---

[1]University of Virginia, Corcoran Department of Philosophy, Charlottesville, Virginia 22904-4780, USA, e-mail: pwh2a@virginia.edu

[2]Although what follows draws on Freedman's views, they should be taken as representing my own position rather than his.

[3]A good account of the development of modern probability theory can be found in Von Plato [25].

Hilbert's sixth problem, but this is incorrect, for in the process of solving part of the problem, Kolmogorov transformed it. Hilbert's own formulation of the sixth problem makes this quite clear: "The investigations on the foundations of geometry suggest the problem: To treat in the same manner, by means of axioms, those *physical sciences* in which mathematics plays an important part; in the first rank are the theory of probabilities and mechanics. As to the axioms of the theory of probabilities, it seems to me desirable that their logical investigation should be accompanied by a rigorous and satisfactory development of the method of mean values in mathematical physics, and in particular in the kinetic theory of gases." Hilbert [10], (emphasis added). This clearly indicates that for Hilbert at least, probability was viewed as a part of science, not of mathematics.[4]

As judges of what counted as mathematics the credentials of Hilbert and Kolmogorov are impeccable, and so probability theory must at some point have made the transition from science to mathematics. The task is to clearly show how empirical content is associated with probability theory while allowing it to retain its status as a part of pure mathematics. One way to reject the view that probability theory itself has empirical content is to view the measure-theoretic formulation of probability theory as a purely formal theory, as a symbolic system that has no interpretation but that imposes formal constraints on the properties of measures, random variables, and other items in the domain of the theory. If one does take this line of response, it makes answers to our second and third questions–How do we get from the formal theory to applications? And what is the correct interpretation of probability theory?–more pressing because one has to inject empirical content that is clearly probabilistic into the formal theory in order to apply it. The relation of probability theory to its applications will be somewhat different from the relation of even quite abstract scientific theories to the world, since the axioms of the latter theories already contain empirical content.[5]

## 3. Formal probability theory

To resolve this problem, and to begin to formulate answers to our three questions, we can turn to a suggestion made in Freedman and Stark [9]. There they claim: "Probability has two aspects. There is a formal mathematical theory, axiomatized by Kolmogorov [15]. And there is an informal theory that connects the mathematics to the world, i.e., defines what 'probability' means when applied to real events." (p.201). This is a promising place to start, although "theory" is a little grand for the second component–something along the lines of "mapping" would be more accurate–and I would demur at the project being about meaning. The project is better construed as one addressing how probability values are correctly assigned within a model. The formal mathematical theory to which they refer is, of course, the theory developed in Kolmogorov's seminal Grundbegriffe der Wahrscheinlichtkeitsrechnung (Kolmogorov [15]) and later developments thereof.[6] This theory is centred on the following familiar apparatus:

---

[4]For a detailed examination of the relations between Kolmogorov's work and Hilbert's sixth problem see Hochkirchen [11]. I thank Michael Stoeltzner for the reference and discussions on the topic.

[5]I shall not here address a famous argument due to Quine [19], the conclusion of which is that no sharp distinction can be drawn between mathematical and scientific theories. A response to that argument involves philosophical issues that will be addressed in a future publication.

[6]There are other, less widely known, examples of this type such as the theories of Alfred Renyi [21], p. 38 and Karl Popper [18], Appendices *iv and *v, both of which take conditional probability rather than absolute probability as a primitive.

Given a set $\Omega$, a $\sigma$-algebra $F$ on $\Omega$ and a real valued set function $P$ on $F$, the triple $<\Omega, F, P>$ constitutes a *probability space* if for any $A_i \in F$,

(i) $P(A_i) \geq 0$.
(ii) If $A_i \bigcap A_j = \emptyset$ for $i \neq j$, then $P(\bigcup_{i=1}^{\infty} A_i) = \sum_{i=1}^{\infty} P(A_i)$.
(iii) $P(\Omega) = 1$.

I say the standard theory is "centred" on this apparatus, because even though many discussions of how to interpret probabilities refer only to these axioms, it would be perversely narrow to identify probability theory with this minimal basis. Kolmogorov, for example, recognized the importance of the role played in his theory by the definitions of stochastic independence and conditional expectations. In addition, from the very beginning of the modern era the theory of stochastic processes has formed an essential part of the theory of probability and supplements such as ergodicity, martingales, exchangeable measures and the like must surely also be included.

As in any formal theory, the intrinsic nature of the elements of the domain of the probability space $\Omega$, and hence of $F$, is irrelevant, a fact that is captured by the use of induced probability spaces. The space is induced by mapping a given probability space $<\Omega, F, P>$ within which the elements of $\Omega$ are actual outcomes (rather than formal representations of them) onto the abstract space $<\mathbb{R}^1, B^1, \mu>$ using a random variable[7] $X$, where $B^1$ is the Euclidean Borel field on $\mathbb{R}^1$ generated by the collection of intervals $(a, b], -\infty < a < b < +\infty$, and $\mu$ is given by:

$$\forall B \in B^1, \mu(B) = P\{X^{-1}(B)\} = P\{\omega \mid X(\omega) \in B\}.$$

The space $<\mathbb{R}^1, B^1, \mu>$ induced by the random variable $X$ is in various generalized versions the canonical object of attention for a good deal of mathematical probability theory. The mathematical advantages provided by the induced probability space lead to the view that random variables simply "re-label" the outcomes of observations and experiments and that all of the essential probabilistic features can be found in the induced space. As a result, in standard presentations of probability theory a sharp separation is not maintained between the representations of the outcomes and the outcomes themselves, a situation which supports the formalist view and the position that if two generating systems ("trials") have the same probability space, then they are probabilistically identical.[8] This is one reason why, considered at this level of abstraction, probability theory can be considered to be a part of pure mathematics – the objects with which it is concerned are mathematical objects that can be taken sui generis.

## 4. Probability models

Returning to Freedman and Stark's article, its immediate purpose was critical. The authors argued that predictions of earthquakes should be viewed with extreme scepticism because the models upon which the probability values were based are imprecise, poorly motivated, and based upon slim empirical evidence. I have no disagreement with that conclusion. There is also a positive message in the article that is muted yet deserves attention. It is this: rather than locating probabilities in the abstract theory or locating probabilities as features of the world, a more realistic approach is to emphasize the role of probabilistic models, and to locate the probabilistic content in those models. One significant advantage of this approach is that it makes much clearer the way in which subject matter specific scientific content plays a role in attributions of probability values.

---

[7] Or in the more general case, sets of random variables.
[8] See, e.g. Itô [13], p.2.

This need to supplement the pure theory with specific models was recognized by Kolmogorov when he wrote in regard to his own theory: "Our system of axioms is not, however, *complete*, for in various problems in the theory of probability different fields of probability have to be examined" Kolmogorov [15], p. 3.[9] What is needed can be seen by noting that probability theory has at least two roles. There is its representational role – specifying what are the general properties possessed by probabilities – and its use as a calculus, allowing the computation of specific probability values. Until some particular measure or distribution is introduced, abstract probability theory has few non-trivial uses as a calculus. In order to keep separate these uses, we can consider the formal theory of probability as a *mathematical template* (see Humphreys [12]) within which specific models and their particular distributions can be substituted. We can usefully draw parallels between the situation in probability theory and that in classical mechanics. In the force-based version of the latter theory, certain fundamental principles are laid down that require the specification of a force function in order for the theory to be applied. Newton's Second Law $F = ma$ places only minimal constraints on forces and has only formal content until a particular force function has been substituted for the place-holder $F$. It is these fundamental principles that are the templates – they are general mathematical forms within which substitution instances can be made for purposes of empirical application. The need for this kind of substitution within highly abstract theories is common – the specification of the form of Lagrangians or Hamiltonians in other versions of classical mechanics, and the specification of basis sets and Hamiltonians in quantum mechanics are but three widely used examples of this need.

In applying the theory to particular systems, some special features of the distributions will ordinarily be used because at the level of the Kolmogorov axioms, probabilities have no internal structure beyond the minimum imposed by the axioms. From the perspective of abstract probability theory this is understandable, for many of the core results in the area are dependent upon the choice of measure or of the probability space only in very general ways, such as requiring the measure to be separable. In contrast, the specification and articulation of particular distributions is of considerable importance for applying probabilities, because the structure of the distribution is often motivated by considerations, however elementary, about the subject matter to which it applied. In what follows, I shall reserve the term "particular distribution" to denote a probability distribution having a density or mass function identifiable by a specific functional form, such as that of the hypergeometric distribution. The plain term "distribution" denotes a more abstract sense in which (in the one dimensional case) a distribution is simply a probability measure over the Borel algebra on the reals. If the mathematical form of the substitution instance is computationally tractable the template becomes a *computational template* and the template can then be used as a calculus. In the example of classical mechanics, the templates form a familiar part of the apparatus of ordinary and partial differential equations, only some of which lend themselves to a closed form solution, given appropriate initial or boundary conditions.

---

[9] By a "field of probability", Kolmogorov meant anything that satisfies the axioms.

## 5. Model and concrete generating systems

Here, again, is Kolmogorov: "We apply the theory of probability to the actual world of experiments in the following manner: (1) There is assumed a complex of conditions, $\mathfrak{C}$, which allows of any number of repetitions. (2) We study a definite set of events which could take place as a result of the establishment of the conditions $\mathfrak{C}$." Kolmogorov [15], p.3. Combining this with the Freedman and Stark approach, we shall need two types of structure in addition to those making up the formal theory in order to capture the relation between probability theory and the world.

First, we have the pair $<MS, MP>$ consisting of the *model generating system* $MS$ and its associated *model probability distribution* $MP$. $MS$ serves as the source of the elements in $\Omega$, the outcome set of the probability space. The structure of $MS$, which can be quite abstract, constrains and sometimes even determines the structure of $MP$ in a way illustrated by the Poisson model described below. $MP$ is the substitution in the probability template that converts the distribution P occurring in the probability space into a particular distribution. Both $MS$ and $MP$ are mathematical objects: "...probability is just a property of a mathematical model intended to describe some features of the natural world... This interpretation – that probability is a property of a mathematical model and has meaning for the world only by analogy – seems the most appropriate for earthquake prediction. To apply the interpretation, one... interprets a number *calculated from the model* to be the probability of an earthquake in some time interval." Freedman and Stark [9], p.5 (emphasis added).

Second, we have the pair $<CS, CP>$ consisting of a *concrete generating system* $CS$ and the corresponding *concrete probability distribution* $CP$. Examples of concrete generating systems are a die thrown under specified conditions, a radioactive atom, and a stock traded on a market. These are real dice, real atoms, and real stock markets, not representations of them. Concrete systems give rise to concrete outcomes but it is very easy to conflate these concrete outcomes with our representations of them. So, to be quite clear: by the outcome of a process such as a die toss, I mean the outcome of a particular side coming up, not a representational description such as "6".

It is useful to think in terms of the hierarchy

$$<\Omega, F, P> \Rightarrow <MS, MP> \Rightarrow <CS, CP>$$

and the intermediate link is the focus of model based probabilities.[10] Everything in this hierarchy except the members of the rightmost element is a formal mathematical object. In virtue of specifying the particular distribution $MP$ and the relation between $MS$ and $\Omega$, one moves from the template at the left to the middle element. The relation between the middle element and the concrete system is discussed with respect to interpretations in Section 8 below.

## 6. Some examples

To see how the hierarchy works, consider the very simple example of a system which has the structure of a Poisson process. This example lies between the simple transparent models of coin tossing and dice throwing and the complex opaque models criticized by Freedman and Stark and so can perhaps better illuminate how model based probability distributions are generated.

---

[10]In many cases there will also be an abstract, non-mathematical, model of $CS$ between $MS$ and $CS$ but for simplicity we can assume that object is used heuristically and is not part of the deductive apparatus.

Here is one set of assumptions behind the attribution of a Poisson process to a system:

(a) During a small interval $t$, the chance of one event occurring in that interval is approximately proportional to the length of that interval: $P(N_t = 1) = \lambda t + f(t)$ where $f(t) \in o(t)$, i.e.

$$\lim_{t \to 0} \frac{P\{N_t = 1\} - \lambda t}{t} = 0.$$

(b) The chance of two or more events occurring in a small interval is small and goes to zero rapidly as the length of the interval goes to zero: $P(N_t \geq 2) \sim o(t)$.

(c) The chance of $n$ events occurring in a given interval is independent of the chance of $m$ events occurring in any disjoint interval, for any $n, m \geq 0$.

(d) The chance of $n$ events occurring in an interval of length $t$ depends only upon the length of the interval and not upon its location.

These facts about the process can be given precise mathematical representations in familiar ways and from those representations one can derive the exact form of the probability distribution covering the output from the process. Within the broad spectrum of applications of the Poisson process, there is a division between those which are justified solely in terms of a reasonable fit between the observed distribution of frequencies and the model distribution, which we can call *frequency driven models*, and those for which some moderately plausible scientific model lies behind the adoption of the probability model, which we can call *theory-driven models*. An example of the former is the analysis of flying bomb hits on London during World War II (Clarke [3]), for which the empirical fit to a Poisson distribution is reasonably good. Yet there is no plausible aerodynamical or military model that would explain why the trajectories of V-1 rockets should satisfy this particular distribution.

In other cases, some attempt is made at providing a theory driven model and one of the appealing features of the Poisson process is that such models can be easily generalized. The successive development of a simple model for fluctuations in electron-photon cascades in absorbers provides an illustration of this. Electron-photon cascades occur when high energy electrons colliding with atoms in an absorber produce photons that lead to production of further electrons, producing a cascade effect. The initial model by Bhabha and Heitler published in 1937 identified $MS$ with a basic Poisson process within which $t$ represents the thickness of the absorber and $n$ represents the number of electrons above an energy level $E$, so that

$$P_n(E, t) = e^{-\lambda t} \frac{(\lambda t)^n}{n!}; \lambda = \lambda(E).$$

However, the predicted mean value for $n$ of $\lambda t$ from this particular model is physically implausible because thicker materials tend to absorb energy from electrons and so it was first modified to a linear birth process and then to a Pólya process. The fact that none of these distributions is completely realistic reflects the severe simplifications resulting from treating the energy levels as discrete and model generating systems of greater sophistication were subsequently developed.[11]

A different example of how a model generating system can be connected with a model probability distribution can be extracted from a result found in Keller [14].

---
[11] See Bharucha-Reid [2].

Keller considered an idealized coin of negligible thickness having its centre of gravity at its geometrical centre (thus making it bias-free). Given the initial conditions of $u =$ upward velocity of the centre of gravity and $\omega =$ angular momentum of a diameter of the coin, Keller showed that when $u \to \infty$ and $\omega \to \infty$ the chance of heads $\to 0.5$, irrespective of what continuous probability density $p(u, \omega)$ describes the initial conditions. $MS$ here consists of the mathematical representation of the idealized coin, the distribution on the initial conditions, the trajectory of the coin, and an absorbing surface. In this case, it would be appropriate to assign to a real coin a chance of coming up heads of approximately 0.5 when large values of $u$ and $\omega$ are present. Note again by what thin threads the middle level mathematical model is tied to the real system, a small amount of idealized physics sufficing to ground the model. Indeed, Diaconis, Holmes, and Montgomery [5] question the applicability of the analysis to real coins, partly on the basis of experimental data, partly by examining the validity of the associated physical model.

## 7. Empirical content

We now have the apparatus to answer our original question: Is probability theory a mathematical theory or a scientific theory? In the three layers of our representation

$$< \Omega, F, P > \Rightarrow < MS, MP > \Rightarrow < CS, CP >$$

the abstract probability space and the accompanying Kolmogorov theory are parts of pure mathematics that have no factual content. The development of Kolmogorov's theory may once have been partially motivated by empirical concerns, and indeed within elementary probability, relative frequencies were a guide for Kolmogorov,[12] but there is no frequency interpretation for the full non-elementary theory that is adequate.[13] In addition, the full theory makes an essential appeal to infinite collections and so has no direct empirical content. Once again, Kolmogorov: "Since the new axiom [Axiom VI, the continuity axiom, from which countable additivity follows] is essential for infinite fields of probability only, it is almost impossible to elucidate its empirical meaning... Infinite fields of probability occur only as idealized models of real random processes. *We limit ourselves, arbitrarily, to only those models which satisfy Axiom VI.*" (Kolmogorov [15], p. 15, italics in original).

Our intermediate element, the model generating system, is also a mathematical object. Consider Freedman and Stark's example (op. cit, p. 205) of a Maxwell–Boltzmann distribution being replaced by a Bose–Einstein distribution for Bose–Einstein condensates. Initially we have a Maxwell–Boltzmann model as the middle element, the properties of which are constrained by the Kolmogorov theory.[14] Now suppose that computer-generated data from a Bose–Einstein distribution are compared with the Maxwell–Boltzmann model. They will fail to fit that model and as a result the Maxwell–Boltzmann model will be replaced by a Bose–Einstein model. Note that everything in this scenario is a mathematical object.[15] In such a case, where the form of the distribution is changed, rather than estimates of particular probabilistic parameters, it is incorrect to say that the Maxwell–Boltzmann model is revised. One constructs or selects a probability distribution for substitution in the template and that distribution either correctly represents the data or it does not. If it does not, the process involves replacement of the model, not its revision. Now consider the situation in which data from a real Bose–Einstein system (e.g. Anderson et al. [1]) are used. Once again, if the predictions from the Maxwell–Boltzmann model do not fit the factual statistics, then one replaces the model with a better particular distribution. That particular distribution has a fixed mathematical form and any changes in its form take us to a different distribution.

---

[12] "In establishing the premises necessary for the applicability of the theory of probability to the world of actual events, the author has used, in large measure, the work of R. von Mises [23], pp. 21–27." Kolmogorov [15], p. 3 (footnote 5.)

[13] See van Fraassen [22], pp. 184–187.

[14] These constraints are usually tacit and the model generating system is often considered to be an autonomous object of investigation.

[15] There is a small factual element if the data are generated on a real machine. We can either ignore that, or think in terms of a virtual machine generating the data.

Here is what we now have: There is the general mathematical template of the Kolmogorov theory and its abstract probability spaces. There is a collection of mathematical models – model generating systems. If you have Platonist inclinations, this collection is very large, it is completely abstract, and it contains models not yet known to us. If your inclinations are more constructivist, the collection will include only models known to us and built by us. All mappings between the abstract probability space and the collection of mathematical models will be mathematical. Finally, there will be a collection of real, concrete generating systems located in the world. There will be mappings between the collection of model generating systems and the collection of concrete generating systems but only a select few will be identified by users and asserted to form the basis of a modeling relation between the mathematical model and the real system. It is these mappings that contain the empirical content, not the mathematical models or the general theory. They are empirical in the sense that empirical facts about the concrete generating system play a role in whether the mapping is structure preserving. The grounds on which these mappings are assessed is by no means simple – the problem of relating physical independence in concrete systems to stochastic independence in the model is by itself a notoriously difficult task – but the point here is that what is false is the assertion that a particular $MS$ and $MP$ have been correctly mapped onto a specific $CS$ and $CP$. The $MS$ and $MP$ involved can then be replaced, but they have no empirical content themselves.

## 8. Interpretations

Finally, what about our third question concerning the interpretation of probability? In the previous section I mentioned the structural features of probability models. What of the specific probability values that are the main concern of the Freedman and Stark paper? Even from the largely formalist perspective adopted here, we cannot ignore the long tradition of trying to provide a substantive interpretation for probability theory for it underlies the differences, sometimes contentious, between subjective Bayesians, objective Bayesians, frequentists, and other schools of thought in probability and statistics. Probability theory is also tied to statistics and whether one chooses to explore the properties of loss functions or to favour classical Neyman-Pearson hypothesis testing, the interpretative issue has to be addressed. And what is perhaps most important, the present approach involves mapping a model probability distribution onto something concrete. What could that something be?

Two different approaches have been used to provide an interpretation for probability theory. The first approach uses an explicit definition of the term "probability" or "has a probability of value $p$". This approach was used by Hans Reichenbach [20], Richard von Mises [24], and Bruno de Finetti [4]. For example, in von Mises' account, as modified by Alonzo Church, we have: Event $A$ has a probability of value $p$ relative to the collective $R$ if and only if $A$ is an event of type $\boldsymbol{A}$, $A$ occurs in $R$, and events of type $\boldsymbol{A}$ occur in $R$ with limiting relative frequency $p$. A collective is an infinite sequence of events within which all event types have limiting relative frequencies that are invariant under selection of subsequences by recursively definable functions on initial segments of the sequence.

Such explicit definitions have the virtue of reducing the concept of probability to other, presumably less opaque, concepts which in the case of the von Mises/Church approach are those of arithmetic limits and recursive functions and it ties the theory based on the definitions very tightly to the intended interpretation. The disadvantage is that such explicit definitions lead to accounts of probability that are different from the account provided by the standard Kolmogorov axiomatization and the theories of which are less general than the measure-theoretic account. For example, de Finetti's theory of personal probabilities, based on the concept of an agent's rational degrees of belief, rejects the countable additivity property of standard probabilities on the grounds that it is operationally meaningless. De Finetti claimed, quite plausibly, that human agents cannot distribute their degrees of belief over infinite sets of outcomes. Von Mises' frequentist theory rejected theorems of standard probability theory about events that occur infinitely often, such as the Borel-Cantelli Lemmas, on the grounds that such theorems were empirically unverifiable.

In contrast, the second approach to interpreting probability uses implicit definitions. Recognizing that chains of definitions must be grounded in primitive terms, this approach takes "probability" as a primitive and relies on a formal theory to place constraints on the probabilities. This second approach has two aspects worth noting. The first is that it captures the idea that all of the specifically probabilistic content is contained in the formal probability spaces. The second aspect is that this approach has the consequence that any probability spaces that are isomorphic are treated as indistinguishable. This is the position underlying the use of the induced probability measures discussed earlier in Section 3. It is for this reason that the abstract theory has only formal content – any attempt to impose a more specific interpretation will be arbitrary.

The term "probability" under this second approach thus refers to an element in a formal theory. At the intermediate level there are particular formal distributions and concrete generating systems can produce statistical estimates of values associated with those distributions. These estimates have an inescapably finite basis and can be generated using finite frequencies, rational degrees of belief, in some cases symmetry arguments, or other means. But these estimates are not interpretations of probability, they are measurements of a parameter's value. The publication of Kolmogorov's *Grundbegriffe* marked a sharp division between the formal theory of probability and those approaches, such as von Mises' and de Finetti's, that used idealizations of methods for estimating the values of elements in models, calling these idealizations "probabilities". The two sides of the division can be brought into only indirect contact.

As an analogy, consider determinism. There are formal theories of determinism that capture this intuition: a system is (historically) deterministic if, under the constraints imposed on that system by the laws that govern it, any complete state of the system is mapped onto a unique later state of the system.[16] There are model systems that are deterministic according to these theories of determinism.[17] And there are concrete systems that, within measurement error, behave similarly to the model systems. So to ask "Are there systems that are deterministic?" is a sensible question with an affirmative answer. But the question "What is determinism?" when asked within this theory is misplaced. All we have is the formal abstract theory together with some specific deterministic systems. And so too with probability. There is the abstract formal theory and there are the various particular probabilistic models. The question "What is probability?" is properly approached through the latter, not the former.

---

[16]See Montague [17].
[17]But considerably fewer than the famous claim in Laplace [16] suggested.

How does factual probabilistic content find its way into the models? There will occur within the model generating system some parameter or distribution representing probabilities, the fact that they are probabilities being grounded in their satisfying the constraints placed on them by the mathematical theory. With frequency driven models, the probabilities will be interpreted as frequencies; with theory driven models, theoretically grounded input will help constrain probability values and distributional forms; perhaps even Bayesian methods can be brought to bear on other types of models. This, I believe, is where the Freedman and Stark model-based probability view constitutes both a distinctive position and a practical but cautionary note. It is distinctive because it directs the interpretative enterprise away from the theory of probability and towards specific probability models. It is practical because it forces one to consider what subject matter specific knowledge is required and is available to inject values into those models. And it is cautionary because it draws our attention to the fact that in many, perhaps most, models the amount of information available is far less than we need to make serious numerical assignments of probability.

**Acknowledgment.** Thanks to David Freedman for his many comments on drafts of this paper. He continues to disagree with some of the claims made here.